\documentclass[a4paper,oneside,notitlepage]{article}%
\usepackage{amsmath}
\usepackage{amssymb}
\usepackage{amscd}
\usepackage{amsfonts}
\usepackage{amstext}
\usepackage{graphicx}
\usepackage{bbm}
\usepackage{bbm}%
\providecommand{\U}[1]{\protect\rule{.1in}{.1in}}
\newtheorem{theorem}{Theorem}

\newtheorem{remark}[theorem]{Remark}

\newtheorem{corollary}[theorem]{Corollary}

\newtheorem{definition}[theorem]{Definition}

\newtheorem{proposition}[theorem]{Proposition}

\begin{document}

\title{Laplace copulas of multifactor gamma distributions are new generalized
Farlie-Gumbel-Morgenstern copulas}
\author{Philippe Bernardoff\\Universit\'{e} de Pau et des Pays de l'Adour\\Laboratoire de Math\'{e}matiques et de leurs Applications, UMR 5142\\avenue de l'Universit\'{e}\\64000 Pau, France.\\e.mail : philippe.bernardoff@univ-pau.fr}
\date{November 10, 2016}
\maketitle

\begin{abstract}
This paper provides bifactor gamma distribution, trivariate gamma distribution
and two copula families on $[0,1]^{n}$ obtained from the Laplace transforms of
the multivariate gamma distribution and the multi-factor gamma distribution
given by $\left[  P\left(  \boldsymbol{\theta}\right)  \right]  ^{-\lambda}$
and \newline$\left[  P\left(  \boldsymbol{\theta}\right)  \right]  ^{-\lambda
}\prod_{i=1}^{n}\left(  1+p_{i}\theta_{i}\right)  ^{-\left(  \lambda
_{i}-\lambda\right)  }$ respectively, where $P$ is an affine polynomial with
respect to the $n$ variables $\theta_{1},\ldots,\theta_{n}$.

These copulas are new generalized Farlie-Gumbel-Morgenstern copulas and allow
in particular to obtain multivariate gamma distributions for which the
cumulative distribution functions and the probability distribution functions
are known.\newline

\rule[0.02cm]{10cm}{0.02cm}


\newpage

\end{abstract}

\section{Introduction}

This paper is motivated by the opportunity to produce explicit cumulative
distribution functions (cdf) and probability distribution functions (pdf) of
multivariate gamma distributions and multi-factor gamma distributions. In this
way, we present two new classes of multivariate copulas generalizing the
Farlie-Gumbel-Morgenstern Copulas. From Sklar (Sklar, 1959) who states that
the cdf $F$ of a random vector $\mathbf{X}=\left(  X_{1},\ldots,X_{n}\right)
$ with continuous marginal cdfs can be uniquely written in the form%
\begin{equation}
F\left(  x_{1},\ldots,x_{n}\right)  =C\left[  F_{1}\left(  x_{1}\right)
,\ldots,F_{n}\left(  x_{n}\right)  \right]  ,\text{ }\mathbf{x}=\left(
x_{1},\ldots,x_{n}\right)  \in\mathbb{R}^{n}, \label{Sklar1}%
\end{equation}
where $n$ is the dimension of the random vector $\mathbf{X}$, where $C:\left[
0,1\right]  ^{n}\rightarrow\left[  0,1\right]  $ is a \emph{copula} and where
$F_{1},\ldots,F_{n}$ are the marginal cdfs of $\mathbf{X}$.

If $f$ is the pdf of the random vector $\mathbf{X}$ and $f_{1},\ldots,f_{n}$
the marginal pdfs of $\mathbf{X,}$ then%
\[
c\left(  u_{1},\ldots,u_{n}\right)  =\frac{\partial^{n}}{\partial u_{1}%
\ldots\partial u_{n}}C\left(  u_{1},\ldots,u_{n}\right)  ,
\]
and we have the following equality%
\begin{equation}
f\left(  x_{1},\ldots,x_{n}\right)  =c\left[  F_{1}\left(  x_{1}\right)
,\ldots,F_{n}\left(  x_{n}\right)  \right]  f_{1}\left(  x_{1}\right)  \ldots
f_{n}\left(  x_{n}\right)  . \label{Sklar2}%
\end{equation}

From Equality (\ref{Sklar1}) it is possible to express the cdf $F$ of
$\mathbf{X}$. This expression cannot simply explicitly give the copula $C$ for
the \textit{multivariate gamma distribution associated with} $\left(
P,\lambda\right)  $ and the \textit{multi-factor gamma distribution associated
with} $\left(  P,\Lambda\right)  .$ But with Joe (Joe, 1997), we can give
copulas deduced of the Laplace transform of the \textit{multivariate gamma
distribution associated with} $\left(  P,\lambda\right)  $ and the Laplace
transform of the \textit{multi-factor gamma distribution associated with}
$\left(  P,\Lambda\right)  .$ By applying the formulas (\ref{Sklar1}) and
(\ref{Sklar2}) we can give a new \textit{multivariate gamma distribution
associated with} $\left(  P,\lambda\right)  $ and a new\textit{\ multi-factor
gamma distribution associated with} $\left(  P,\Lambda\right)  $ for which we
have an explicit formula for its cdf and its pdf.

The paper is organized as follows. Section 2 gives definitions of multivariate
gamma distributions and multi-factor gamma distributions for which Laplace
transform is given, and considers the bidimentional and tridimentional case.
Section 3 defines the Laplace copula. Section 4 states the two main results.
Section 5 applies the main cases to bidimensional and tridimensional cases. In
particular, the Joe's family BB10 is generalized. For ease of the fluent
exposition of the paper, proofs are collected in the Appendix.

\section{Multivariate gamma distributions and multi-factor\textit{\ }gamma
distributions}

In the literature, the multivariate gamma distributions on $\mathbb{R}^{n}$
have several non-equivalent definitions. Many authors require only that the
marginal distributions are ordinary gamma distributions (Balakrishnan
\textit{et al}., 1997). In the present paper, we use the extension of the
classical one-dimensional definition to $\mathbb{R}^{n}$ obtained as follows:
we consider an affine polynomial $P\left(  \boldsymbol{\theta}\right)  $ in
$\boldsymbol{\theta}=\left(  \theta_{1},\ldots,\theta_{n}\right)  $ where
`affine' means that, for $j=1,\ldots,n,$ $\partial^{2}P/\partial\theta_{j}%
^{2}=0$. We also assume that $P\left(  \mathbf{0}\right)  =1.$ For instance,
for $n=2,$ we have $P\left(  \theta_{1},\theta_{2}\right)  =1+p_{1}\theta
_{1}+p_{2}\theta_{2}+p_{\left\{  1,2\right\}  }\theta_{1}\theta_{2}$.

We denote by $\mathfrak{P}_{n}=\mathfrak{P}\left(  \left[  n\right]  \right)
$ the family of all subsets of $\left[  n\right]  $ and $\mathfrak{P}%
_{n}^{\ast} $ the family of non-empty subsets of $\left[  n\right]  =\left\{
1,\ldots,n\right\}  .$ For simplicity, if $n$ is fixed and if there is no
ambiguity, we denote these families by $\mathfrak{P}$ and $\mathfrak{P}^{\ast
}$, respectively.

We denote by $\mathbb{N}$ the set of non-negative integers. If $\mathbf{z}%
=\left(  z_{1},\ldots,z_{n}\right)  \in\mathbb{R}^{n}$ and $\mathbf{\alpha
}=\left(  \alpha_{1},\ldots,\alpha_{n}\right)  \in\mathbb{N}^{n},$ then
$\mathbf{\alpha}!=\alpha_{1}!\ldots\alpha_{n}!$, $\left\vert \mathbf{\alpha
}\right\vert =\alpha_{1}+\ldots+\alpha_{n}$, $a_{\mathbf{\alpha}}%
=a_{\alpha_{1},\ldots,\alpha_{n}}$ and
\[
\mathbf{z}^{\mathbf{\alpha}}=\prod\limits_{i=1}^{n}z_{i}^{\alpha_{i}}%
=z_{1}^{\alpha_{1}}\ldots z_{n}^{\alpha_{n}}.
\]
For $T\ $in $\mathfrak{P}_{n},$ we simplify the above notation by writing
$\mathbf{z}^{T}=\prod\nolimits_{t\in T}z_{t}$ instead of $\mathbf{z}%
^{\mathbf{1}_{T}}$ where
\[
\mathbf{1}_{T}=\left(  \alpha_{1},\ldots,\alpha_{n}\right)  \text{ with
}\alpha_{i}=1\text{ if }i\in T\text{ and }\alpha_{i}=0\text{ if }i\notin T.
\]
We also write $\mathbf{z}^{-T}$ for $\prod\nolimits_{t\in T}1/z_{t}.$ For a
mapping $a:\mathfrak{P}\rightarrow\mathbb{R}$, we shall use the notation
$a:\mathfrak{P}\rightarrow\mathbb{R}$, $T\mapsto a_{T}.$ In this notation, an
affine polynomial with constant term equal to $1$ is $P\left(
\boldsymbol{\theta}\right)  =\sum_{T\in\mathfrak{P}}p_{T}\boldsymbol{\theta
}^{T},$ with $p_{\varnothing}=1.$ For simplicity, if $T=\left\{  t_{1}%
,\ldots,t_{k}\right\}  ,$ we denote $a_{\left\{  t_{1},\ldots,t_{k}\right\}
}=a_{t_{1}\ldots t_{k}}.$ The indicator function of a set $S$ is denoted by
$\mathbbm{1}_{S},$ that is, $\mathbbm{1}_{S}\left(  x\right)  =1$ for $x\in S$
and $0$ for $x\notin S$.

We fix $\lambda>0.$ If a random vector $\mathbf{X=}\left(  X_{1},\ldots
,X_{n}\right)  $ on $\mathbb{R}^{n}$ with probability distribution (pd)
$\mu_{\mathbf{X}}$ is such that its Laplace transform is%
\begin{equation}
\mathbb{E}\left\{  \exp\left[  -\left(  \theta_{1}X_{1}+\cdots+\theta_{n}%
X_{n}\right)  \right]  \right\}  =\left[  P\left(  \boldsymbol{\theta}\right)
\right]  ^{-\lambda}, \label{TLMGD}%
\end{equation}
where $\mathbb{E}$ denotes the expectation, for a set of $\boldsymbol{\theta}$
with non-empty interior, then we denote $\mu_{\mathbf{X}}=\gamma_{\left(
P,\lambda\right)  },$ and $\gamma_{\left(  P,\lambda\right)  }$ will be called
the \textit{multivariate gamma distribution associated with} $\left(
P,\lambda\right)  $. These multivariate gamma distributions occur naturally in
the classification of natural exponential families in $\mathbb{R}^{n}$
(Bar-Lev \textit{et al.,} 1994).

The marginal distributions of the \textit{multivariate gamma distribution
associated with} $\left(  P,\lambda\right)  $ are ordinary gamma distributions
of parameters $\left(  p_{i},\lambda\right)  $ for $i=1,\ldots,n,$ with
Laplace transform $\left(  1+p_{i}\theta_{i}\right)  ^{-\lambda},$ and pd
$\gamma_{\left(  p_{i},\lambda\right)  }\left(  \text{\emph{d}}x\right)
=x^{\lambda-1}p_{i}^{-\lambda}/\Gamma\left(  \lambda\right)  \exp\left(
-x/p_{i}\right)  \mathbbm{1}_{\left(  0,\infty\right)  }\left(  x\right)
$\emph{d}$x.$ We extend the first definition to the\textit{\ multi-factor
gamma distribution associated with} $\left(  P,\Lambda\right)  $ where
$\Lambda=\left(  \lambda,\lambda_{1},\ldots,\lambda_{n}\right)  $ and
$\lambda_{i}\geqslant\lambda>0$ for all $i=1,\ldots,n$ by its Laplace
transform%
\begin{equation}
\mathbb{E}\left\{  \exp\left[  -\left(  \theta_{1}X_{1}+\cdots+\theta_{n}%
X_{n}\right)  \right]  \right\}  =\left[  P\left(  \boldsymbol{\theta}\right)
\right]  ^{-\lambda}\prod_{i=1}^{n}\left(  1+p_{i}\theta_{i}\right)
^{-\left(  \lambda_{i}-\lambda\right)  }. \label{TLMMGD}%
\end{equation}
The marginal distributions of the \textit{multi-factor gamma distribution
associated with} $\left(  P,\Lambda\right)  $ are ordinary gamma distributions
of parameters $\left(  p_{i},\lambda_{i}\right)  $ for $i=1,\ldots,n.$

We state first a proposition, whose proof is obvious.

\begin{proposition}
\label{buildMFGD}A random vector $\mathbf{X}$ with distribution
$\mathbf{\gamma}_{P,\Lambda}$ can be obtained in the following way: \newline
Let $\mathbf{Y}$ be a vector with distribution $\mathbf{\gamma}_{P,\lambda}.$
Let $\mathbf{Z}=\left(  Z_{1},\dots,Z_{n}\right)  $ be a random vector
constituted of independent margins for which its pds are $\gamma_{\left(
p_{i},\lambda_{i}\right)  },$ and such that $Z$ and $Y$ are independent. Then
the vector $\mathbf{X}=\mathbf{Y}+\mathbf{Z}$ has Laplace transform
(\ref{TLMMGD}), and consequently is a\textit{\ multi-factor gamma distribution
}associated to $\left(  P,\Lambda\right)  .$
\end{proposition}

For the bidimensional case, Dussauchoy and Berland, (Dussauchoy and Berland,
1972) consider the random vector $\mathbf{X}=(X_{1},X_{2})$ with Laplace
transform
\begin{align}
L_{\mathbf{X}}\left(  \theta_{1},\theta_{2}\right)   &  =\left(  1+p_{1}%
\theta_{1}\right)  ^{-\lambda_{1}}\left(  1+p_{2}\theta_{2}\right)
^{-\lambda_{2}}\left[  1-\frac{r_{12}p_{1}p_{2}\theta_{1}}{\left(
1+p_{1}\theta_{1}\right)  \left(  1+p_{2}\theta_{2}\right)  }\right]
^{-\lambda}\nonumber\\
&  =\left(  1+p_{1}\theta_{1}\right)  ^{-\left(  \lambda_{1}-\lambda\right)
}\left(  1+p_{2}\theta_{2}\right)  ^{-\left(  \lambda_{2}-\lambda\right)
}\left(  1+p_{1}\theta_{1}+p_{2}\theta_{2}+p_{12}\theta_{1}\right)
^{-\lambda}, \label{DuchaussoyBerland}%
\end{align}
where $r_{12}=1-p_{12}/\left(  p_{1}p_{2}\right)  >0$ and $p_{1},p_{2}%
,p_{12}>0 $.

Unfortunately the Laplace transform of these pds is simple, but its pdfs and
cdfs are unknown, except for the case $n=2$ for the \textit{multivariate gamma
distribution associated with} $\left(  P,\lambda\right)  $. Let $F_{m}^{p}$ be
the generalized hypergeometric function (Slater, 1966) defined by%
\begin{equation}
F_{m}^{p}\left(  \alpha_{1},\ldots,\alpha_{p};\beta_{1},\ldots,\beta
_{m};z\right)  =\sum_{k=0}^{\infty}\frac{\left(  \alpha_{1}\right)  _{k}%
\cdots\left(  \alpha_{p}\right)  _{k}}{\left(  \beta_{1}\right)  _{k}%
\cdots\left(  \beta_{m}\right)  _{k}}\frac{z^{k}}{k!}, \label{Hypergeomconf}%
\end{equation}
where $\left(  a\right)  _{k}=\Gamma\left(  a+k\right)  /\Gamma\left(
a\right)  $ for $a>0$ and $k\in\mathbb{N}$, or more generally by $(a)_{0}=1$,
$(a)_{n+1}=(a+n)(a),$ $\forall n\in\mathbb{N}$, $\forall a\in\mathbb{R}$, is
the Pochhammer's symbol. For simplification, we denote $F_{m}^{0}$ by $F_{m}$.
Bernardoff (Bernardoff, 2006) gives the

\begin{proposition}
Let $P\left(  \theta_{1},\theta_{2}\right)  =1+p_{1}\theta_{1}+p_{2}\theta
_{2}+p_{12}\theta_{1}\theta_{2}$ be an affine polynomial where $p_{1},p_{2}>0$
and $p_{1,2}>0.$ Let $\mu=\boldsymbol{\gamma}_{\left(  P,\lambda\right)  }$ be
the gamma distribution associated to $\left(  P,\lambda\right)  .$ The measure
$\mu$ exists if and only if $c=\left(  p_{1}p_{2}-p_{12}\right)  /p_{12}%
^{2}=p_{1}p_{2}/p_{12}^{2}r_{12}>0$. Then we have%
\begin{equation}
\boldsymbol{\gamma}_{\left(  P,\lambda\right)  }\left(  \text{\emph{d}}%
x_{1},\text{\emph{d}}x_{2}\right)  =\frac{p_{12}^{-\lambda}}{\Gamma\left(
\lambda\right)  ^{2}}\mathbf{e}^{-\frac{p_{2}}{p_{12}}x_{1}-\frac{p_{1}%
}{p_{12}}x_{2}}\left(  x_{1}x_{2}\right)  ^{\lambda-1}F_{1}\left(
\lambda;cx_{1}x_{2}\right)  \mathbbm{1}_{\left(  0,\infty\right)  ^{2}}\left(
\mathbf{x}\right)  \,\text{\emph{d}}\left(  \mathbf{x}\right)  . \label{mud2}%
\end{equation}

\end{proposition}

For the case $\Lambda=\left(  \lambda,\lambda,\lambda_{2}\right)  $, the
\textit{multi-factor gamma distribution associated with} $\left(
P,\Lambda\right)  $ is named by Chatelain \emph{et al.} (Chatelain \emph{et
al., }2008) the\textit{\ multisensor gamma distribution associated with}
$\left(  P,\lambda,\lambda_{2}\right)  $ and they have proved that its pd is
given by the equality%
\begin{align}
\boldsymbol{\gamma}_{\left(  P,\Lambda\right)  }\left(  \text{d}x_{1}%
,\text{d}x_{2}\right)   &  =\frac{p_{12}^{-\lambda}p_{2}^{-\left(  \lambda
_{2}-\lambda\right)  }1}{\Gamma\left(  \lambda\right)  \Gamma\left(
\lambda_{2}\right)  }x_{1}^{\lambda-1}x_{2}^{\lambda_{2}-1}\mathbf{e}%
^{-\frac{p_{2}}{p_{12}}x_{1}-\frac{p_{1}}{p_{12}}x_{2}}\Phi_{3}(\lambda
_{2}-\lambda;\lambda_{2};c\frac{p_{12}}{p_{2}}x_{2};cx_{1}x_{2})\nonumber\\
&  \times\mathbbm{1}_{\left(  0,\infty\right)  ^{2}}\left(  x_{1}%
,x_{2}\right)  \,\text{d}x_{1}\text{d}x_{2}\text{,} \label{BMGD}%
\end{align}
where%
\begin{equation}
\Phi_{3}\left(  a;b;x,y\right)  =\sum_{m,n\geqslant0}\frac{\left(  a\right)
_{m}}{\left(  b\right)  _{m+n}}\frac{x^{m}}{m!}\frac{y^{n}}{n!}
\label{HornFunction}%
\end{equation}
is the Horn function.

For the bidimensional general case, we have the following Theorem. Let $F_{I}$
be the function defined by%
\begin{equation}
F_{I}\left(  a,b,c,\mathbf{z}\right)  =\sum_{m_{1},m_{2},m_{3}=0}^{\infty
}\frac{\left(  a\right)  _{m_{1}}\left(  b\right)  _{m_{2}}\left(  c\right)
_{m_{3}}}{\left(  a+c\right)  _{m_{1}+m_{3}}\left(  b+c\right)  _{m_{2}+m_{3}%
}}\frac{z_{1}^{m_{1}}}{m_{1}!}\frac{z_{2}^{m_{2}}}{m_{2}!}\frac{z_{3}^{m_{3}}%
}{m_{3}!}; \label{genLauricella1}%
\end{equation}
it is a particular generalized Lauricella function defined, by example, in
(Panda, 1973).

\begin{theorem}
\label{Proposition_multifactorbig2}The pd of $\boldsymbol{\gamma}_{\left(
P,\left(  \lambda,\lambda_{1},\lambda_{2}\right)  \right)  }$ is given by the
equality%
\begin{align}
\boldsymbol{\gamma}_{\left(  P,\left(  \lambda,\lambda_{1},\lambda_{2}\right)
\right)  }\left(  \text{\emph{d}}x_{1},\text{\emph{d}}x_{2}\right)   &
=\frac{p_{12}^{-\lambda}p_{1}^{-\left(  \lambda_{1}-\lambda\right)  }%
p_{2}^{-\left(  \lambda_{2}-\lambda\right)  }}{\Gamma\left(  \lambda
_{1}\right)  \Gamma\left(  \lambda_{2}\right)  }x_{1}^{\lambda_{1}-1}%
x_{2}^{\lambda_{2}-1}\mathbf{e}^{-\left(  \frac{p_{2}}{p_{12}}x_{1}%
+\frac{p_{1}}{p_{12}}x_{2}\right)  }\nonumber\\
&  \hspace*{-3cm}\times F_{I}\left(  \lambda_{1}-\lambda,\lambda_{2}%
-\lambda,\lambda,\frac{p_{12}}{p_{1}}x_{1},\frac{p_{12}}{p_{2}}x_{2}%
,cx_{1}x_{2}\right)  \mathbbm{1}_{\left(  0,\infty\right)  ^{2}}\left(
x_{1},x_{2}\right)  \,\text{\emph{d}}x_{1}\text{\emph{d}}x_{2}\text{,}
\label{mulbigamma2}%
\end{align}

\end{theorem}

If we get $\lambda_{1}=\lambda$ in the equality (\ref{mulbigamma2}), we obtain
Chatelain and Tourneret's result (\ref{BMGD}) because%
\begin{align*}
F_{I}\left(  0,\lambda_{2}-\lambda,\lambda,z_{1},z_{2},z_{3}\right)   &  =\\
\hspace*{-3cm}\sum_{m_{2},m_{3}=0}^{\infty}\frac{\left(  b\right)  _{m_{2}}%
}{\left(  b+c\right)  _{m_{2}+m_{3}}}\frac{z_{2}^{m_{2}}}{m_{2}!}\frac
{z_{3}^{m_{3}}}{m_{3}!}  &  =\Phi_{3}\left(  b;b+c;z_{2},z_{3}\right)  .
\end{align*}

Bernardoff (Bernardoff, 2006) gives the following Proposition:

\begin{proposition}
Let $\mu$ be a \textit{multivariate gamma distribution }on $\mathbb{R}^{n}$
associated to $\left(  P,\lambda\right)  .$ Assume that $\mu$ is not
concentrated on a linear subspace of $\mathbb{R}^{n}$ of the form
$\{x\in\mathbb{R}^{n};\ x_{k}=0\}$ for some $k$ in $\left[  n\right]
=\{1,\ldots,n\}.$ Then:

\begin{enumerate}
\item[\emph{(i)}] For all $i\in\left[  n\right]  ,$ $p_{i}\neq0$.

\item[\emph{(ii)}] If $p_{1},\ldots,p_{k}<0$ and $p_{k+1},\ldots,p_{n}>0,$
then $\text{Supp}\left(  \mu\right)  \subset\left(  -\infty,0\right]
^{k}\times\left[  0,\infty\right)  ^{n-k}$.

\item[\emph{(iii)}] If $p_{1},\ldots,p_{n}>0$ then $p_{\left[  n\right]
}\geqslant0.$
\end{enumerate}
\end{proposition}

Bernardoff (Bernardoff, 2006) gives a necessary and sufficient condition for
infinite divisibility of the \textit{multivariate gamma distribution
associated with} $\left(  P,\lambda\right)  $, in the sense that the Laplace
transform of $\gamma_{\left(  P,\lambda\right)  }$ power $t$ for all positive
$t$ is still the Laplace transform of a positive measure, by the following theorem:

\begin{theorem}
\label{mainresult} Let $\mu=\gamma_{P,\lambda}$ be a gamma distribution
associated with $\left(  P,\lambda\right)  ,$ where $\lambda>0$ and $P\left(
\boldsymbol{\theta}\right)  =\sum\nolimits_{T\in\mathfrak{P}_{n}}%
p_{T}\boldsymbol{\theta}^{T}$ is such that $p_{i}>0$, for all $i\in\left[
n\right]  ,$ and $p_{\left[  n\right]  }>0$. Let $\widetilde{P}\left(
\boldsymbol{\theta}\right)  =\sum\nolimits_{T\in\mathfrak{P}_{n}}%
\widetilde{p}_{T}\boldsymbol{\theta}^{T}$ be the affine polynomial such that
$\widetilde{p}_{T}=-p_{\overline{T}}/p_{\left[  n\right]  }$ for all
$T\in\mathfrak{P}_{n},$ where $\overline{T}=\left[  n\right]  \smallsetminus
T$. Let
\[
\widetilde{b}_{S}=b_{S}(\widetilde{P})=\sum_{k=1}^{\left\vert S\right\vert
}\left(  k-1\right)  !\sum_{\mathcal{T\in}\Pi_{S}^{k}}\prod_{T\in\mathcal{T}%
}\widetilde{p}_{T},
\]
where $\left\vert S\right\vert $ denotes the cardinality of the set $S$, and
$\Pi_{S}^{k}$ denotes the set of all partitions of $S$ into k \emph{non-empty}
subsets of $S$. Then the measure $\mu$ is infinitely divisible if and only if%
\begin{equation}
\widetilde{\text{$p$}}_{i}<0\text{ for all }i\in\left[  n\right]  ,
\label{cond1}%
\end{equation}
and%
\begin{equation}
\widetilde{b}_{S}\geqslant0\text{ for all }S\in\mathfrak{P}_{n}^{\ast}\text{
such that }\left\vert S\right\vert \geqslant2. \label{cond2}%
\end{equation}

\end{theorem}

\begin{corollary}
By the properties of infinite divisible distributions we conclude that the
necessary and sufficient conditions for infinite divisibility of a
\textit{multivariate gamma distribution }associated to $\left(  P,\lambda
\right)  $ of theorem (\ref{mainresult}), are also necessary and sufficient
conditions for infinite divisibility of \textit{multivariate multi-factor
gamma distribution }associated to $\left(  P,\Lambda\right)  .$
\end{corollary}

To illustrate the difficulty to calculate the \textit{multivariate gamma
distribution associated with} $\left(  P,\lambda\right)  $ we give, for the
tridimensional case, the following theorem. Let $F_{II}$ be the function
defined by%
\begin{equation}
F_{II}\left(  \lambda_{1},\lambda_{2},z_{1},z_{2},z_{3},z_{4}\right)
=\sum_{m_{1},\ldots,m_{4}=0}^{\infty}\frac{1}{\left(  \lambda_{1}\right)
_{m_{1}+m_{2}+m_{3}}\left(  \lambda_{2}\right)  _{2m_{1}+m_{2}+m_{4}}}%
\frac{z_{1}^{m_{1}}}{m_{1}!}\frac{z_{2}^{m_{2}}}{m_{2}!}\frac{z_{3}^{m_{3}}%
}{m_{3}!}\frac{z_{4}^{m_{4}}}{m_{4}!}; \label{genLauricella2}%
\end{equation}
it is still a particular generalized Lauricella function.

\begin{theorem}
\label{gamma3}In the case $n=3$, $p_{i}>0$ for $i\in\left[  3\right]  ,$
$p_{ij}>0$ for $\left(  i,j\right)  \in\left[  3\right]  ^{2}$, $\widetilde{b}%
_{ij}=-\frac{b_{k}}{p_{123}}+\frac{p_{jk}p_{ik}}{p_{123}^{2}}\geqslant0$ for
$i\neq j$ and $\left\{  i,j,k\right\}  =\left[  3\right]  ,$ $p_{123}>0,$ and
\[
\widetilde{b}_{123}=-\frac{1}{p_{123}}+\frac{p_{12}p_{1}}{p_{123}^{2}}%
+\frac{p_{13}p_{2}}{p_{123}^{2}}+\frac{p_{23}p_{1}}{p_{123}^{2}}+2\frac
{p_{12}p_{13}p_{23}}{p_{123}^{3}}\geqslant0,
\]
the infinitely divisible \textit{multivariate gamma distribution }%
$\gamma_{\left(  \lambda,P\right)  }$ \textit{associated with} $\left(
P,\lambda\right)  ,$ is given by the formula%
\begin{align}
\gamma_{\left(  \lambda,P\right)  }\left(  \text{\emph{d}}\mathbf{x}\right)
&  =\frac{p_{123}^{-\lambda}}{\left[  \Gamma\left(  \lambda\right)  \right]
^{3}}\exp(\widetilde{p}_{1}x_{1}+\widetilde{p}_{2}x_{2}+\widetilde{p}_{3}%
x_{3})\left(  x_{1}x_{2}x_{3}\right)  ^{\lambda-1}\nonumber\\
&  \times F_{II}(\lambda,\lambda,\widetilde{b}_{13}x_{1}x_{3}\widetilde{b}%
_{23}x_{2}x_{3},\widetilde{b}_{123}x_{1}x_{2}x_{3},\widetilde{b}_{12}%
x_{1}x_{2},\widetilde{b}_{13}x_{1}x_{3}+\widetilde{b}_{23}x_{2}x_{3}%
)\mathbbm{1}_{\left(  0,\infty\right)  ^{3}}\left(  \mathbf{x}\right)
\,\text{\emph{d}}\mathbf{x}. \label{MGD3}%
\end{align}

\end{theorem}

\begin{remark}
The case $p_{123}=0$ is solved by Letac and Weso\l owski (Letac and
Weso\l owski, 2008).
\end{remark}

\begin{remark}
If $\widetilde{b}_{12}=\widetilde{b}_{13}=\widetilde{b}_{23}=0,$ Theorem
\ref{gamma3} gives
\[
\gamma_{\left(  \lambda,P\right)  }\left(  \text{\emph{d}}\mathbf{x}\right)
=\frac{p_{123}^{-\lambda}}{\left[  \Gamma\left(  \lambda\right)  \right]
^{3}}\exp(\widetilde{p}_{1}x_{1}+\widetilde{p}_{2}x_{2}+\widetilde{p}_{3}%
x_{3})\left(  x_{1}x_{2}x_{3}\right)  ^{\lambda-1}F_{2}\left(  \lambda
,\lambda;\widetilde{b}_{123}x_{1}x_{2}x_{3}\right)  \mathbbm{1}_{\left(
0,\infty\right)  ^{3}}\left(  \mathbf{x}\right)  \,\text{\emph{d}}\mathbf{x},
\]
and if we put $\lambda=1$ in this last equality, we obtain the Kibble and
Moran distribution given in (Balakrishnan \emph{et al.}, 2000).
\end{remark}

\section{Laplace copula}

We recall the following theorem (Marshall and Olkin, 1988), Theorem 2.1, p.
835, and his Corollary

\begin{theorem}
\label{MarshallO}Let $H_{1},\ldots,H_{n},$ be univariate cdfs, and let $G$ be
an $n$-variate cdf such that $\overline{G}\left(  0,\ldots,0\right)  =1,$ with
univariate margins $G_{i}$ $(i=1,\ldots,n)$. Denote the Laplace transforms of
$G$ and $G_{i}$, respectively, by $\phi$ and $\phi_{i}$ $(i=1,\ldots,n)$. Let
$K$ be a copula. If $F_{i}\left(  x\right)  =\exp\left\{  -\phi_{i}%
^{-1}\left[  H_{i}\left(  x\right)  \right]  \right\}  $ $\left(
i=1,\ldots,n\right)  ,$ then
\begin{equation}
H\left(  x_{1},\ldots,x_{n}\right)  =\int\cdots\int K\left\{  \left[
F_{1}\left(  x_{1}\right)  \right]  ^{\theta_{1}},\ldots,\left[  F_{n}\left(
x_{n}\right)  \right]  ^{\theta_{n}}\right\}  \emph{d}G\left(  \theta
_{1},\ldots,\theta_{n}\right)  \label{Formula 1 Marshall and Olkin}%
\end{equation}
is an $n$-variate cdf with marginals $H_{1},\ldots,H_{n}.$
\end{theorem}

\begin{corollary}
\label{corollaireMO1}Under the conditions of the Theorem \ref{MarshallO},
with\newline$K\left(  x_{1},\ldots,x_{n}\right)  =\prod_{i=1}^{n}x_{i},\left(
0\leqslant x_{i}\leqslant1;\text{ }i=1,\ldots,n\right)  ,$%
\begin{equation}
H\left(  x_{1},\ldots,x_{n}\right)  =\phi\left\{  \phi_{1}^{-1}\left[
H_{1}\left(  x_{1}\right)  \right]  ,\ldots,\phi_{n}^{-1}\left[  H_{n}\left(
x_{n}\right)  \right]  \right\}  \label{Formula 2 Marshall and Olkin}%
\end{equation}
defines an $n$-variate cdf with marginals $H_{1},\ldots,H_{n}.$
\end{corollary}

Therefore copula for $H$ is associated to $G$ and is given by the formula%
\begin{equation}
\mathbf{C}\left(  u_{1},\ldots,u_{n}\right)  =\phi\left[  \phi_{1}^{-1}\left(
u_{1}\right)  ,\ldots\phi_{n}^{-1}\left(  u_{n}\right)  \right]  .
\label{Copula Marshall and Olkin}%
\end{equation}

\begin{remark}
If we know $\phi$ and $\phi_{i}$, for $i=1,\ldots,n,$ we can find the copula
$\mathbf{C,}$ that is the case for the multivariate gamma distributions and
for the multivariate multi-factor gamma distributions. For the case
$H_{i}=G_{i},$ we can write :
\begin{equation}
H\left(  x_{1},\ldots,x_{n}\right)  =\phi\left\{  \phi_{1}^{-1}\left[
G_{1}\left(  x_{1}\right)  \right]  ,\ldots\phi_{n}^{-1}\left[  G_{n}\left(
x_{n}\right)  \right]  \right\}  , \label{other gamma}%
\end{equation}
and this last formula defines $n$-variate cdf gamma $H$ with margin
distributions $G_{1},\ldots,G_{n}.$\newline Then the copula associated to $H$
is
\[
\mathbf{C}\left(  u_{1},\ldots,u_{n}\right)  =\phi\left[  \phi_{1}^{-1}\left(
u_{1}\right)  ,\ldots,\phi_{n}^{-1}\left(  u_{n}\right)  \right]  .
\]

\end{remark}

Now, from the Formula \ref{Copula Marshall and Olkin}, we can give the
following Definition

\begin{definition}
Let $\mathbf{X}=\left(  X_{1},X_{2},\ldots,X_{n}\right)  $ be a random vector
in $\left]  0,+\infty\right[  ^{n}$ with Laplace transform $\varphi
_{\mathbf{X}}$ defined by $L_{\mathbf{X}}\left(  \boldsymbol{\theta}\right)
=L_{\mathbf{X}}\left(  \theta_{1},\ldots,\theta_{n}\right)  =\varphi
_{\mathbf{X}}\left(  \theta_{1},\ldots,\theta_{n}\right)  =\mathbb{E}\left\{
\exp\left[  -\left(  \theta_{1}X_{1}+\cdots+\theta_{n}X_{n}\right)  \right]
\right\}  $, and let $\varphi_{X_{i}}$ be the Laplace transform of the random
variable $X_{i},$ defined by $L_{X_{i}}\left(  \theta_{i}\right)
=\varphi_{X_{i}}\left(  \theta_{i}\right)  =\mathbb{E}\left[  \exp\left(
-\theta_{i}X_{i}\right)  \right]  .$ The function $\varphi_{X_{i}}$ is a one
and onto, and is decreasing of $\left[  0,+\infty\right[  $ onto $\left]
0,1\right]  $. Its inverse function is denoted by $\varphi_{X_{i}}^{-1}$. Then
the function $C_{L_{_{\mathbf{X}}}}$ defined by%
\[
C_{L_{_{\mathbf{X}}}}\left(  u_{1},\ldots,u_{n}\right)  =\varphi_{\mathbf{X}%
}\left[  \varphi_{\mathbf{X}_{1}}^{-1}\left(  u_{1}\right)  ,\ldots
,\varphi_{\mathbf{X}_{n}}^{-1}\left(  u_{n}\right)  \right]
\]
is a copula. We call this copula $C_{L_{_{\mathbf{X}}}}$ the Laplace copula
associated to the random vector $\mathbf{X}$. If $\mathbf{X}$ has pd
$\mu_{\mathbf{X}}\left(  \text{\emph{d}}\mathbf{x}\right)  ,$ then we denote
still $C_{L_{_{\mathbf{X}}}}\left(  u_{1},\ldots,u_{n}\right)
=C_{L_{_{\mathbf{\mu}_{\mathbf{X}}}}}\left(  u_{1},\ldots,u_{n}\right)  .$
\end{definition}

\begin{proposition}
Let $\mathbf{F}_{1},\ldots,\mathbf{F}_{n},$ be $n$ univariate cdfs on $\left]
0,+\infty\right[  .$ The relation%
\begin{align*}
\mathbf{F}\left(  x_{1},\ldots,x_{n}\right)   &  =\varphi_{\mathbf{X}}\left\{
\varphi_{\mathbf{X}_{1}}^{-1}\left[  \mathbf{F}_{1}\left(  x_{1}\right)
\right]  ,\ldots,\varphi_{\mathbf{X}_{n}}^{-1}\left[  \mathbf{F}_{n}\left(
x_{n}\right)  \right]  \right\} \\
&  =C_{L_{_{\mathbf{X}}}}\left[  \mathbf{F}_{1}\left(  x_{1}\right)
,\ldots,\mathbf{F}_{n}\left(  x_{n}\right)  \right]
\end{align*}
defines a cdf $\mathbf{F}$ with marginal cdfs $\mathbf{F}_{1},\ldots
,\mathbf{F}_{n}.$ Then the copula associated to $\mathbf{F}$ is
$C_{L_{_{\mathbf{X}}}}.$\newline If we choose, for $i=1,2,\ldots,n$,
$\mathbf{F}_{i}=\mathbf{F}_{X_{i}},$ the cdf of $X_{i},$ then $\mathbf{F}$ is
the cdf of a random vector $\mathbf{X}$ for which the marginal cdfs are the
cdfs of $X_{i}$, and we have $\mathbf{F}\neq\mathbf{F}_{\mathbf{X}}$ where
$\mathbf{F}_{\mathbf{X}}$ calls the cdf of $\mathbf{X}$. In the case of
multivariate gamma distributions and multi-factor gamma distributions, we get
other multivariate gamma distributions and multi-factor gamma distributions
with given copula.
\end{proposition}

Laplace copula can be seen as a generalization of Archimedean copula (Joe, 1997,2014).

\section{Main results}

Now, we are in the capacity to give the two main results of this paper

\begin{theorem}
\label{copulagammaM}Let $P$ an affine polynomial in the $n$ variables
$\theta_{i},$ $i=1,\ldots,n$, with $P\left(  \mathbf{0}\right)  =1,$ then the
Laplace copula of the multivariate gamma distribution $\gamma_{\left(
P,\lambda\right)  }$ associated to $\left(  P,\lambda\right)  $ is%
\begin{equation}
C_{L_{\gamma_{\left(  P,\lambda\right)  }}}\left(  \mathbf{v}\right)
=\mathbf{v}^{\left[  n\right]  }[1+\sum_{T\subset\left[  n\right]  ,\left\vert
T\right\vert >1}\left(  -1\right)  ^{\left\vert T\right\vert }P(-\frac
{1}{\mathbf{p}}\mathbf{1}_{T})\prod_{t\in T}(1-v_{t}^{\frac{1}{\lambda}%
})]^{-\lambda}, \label{LC_MGD}%
\end{equation}
where $\mathbf{v}=\left(  v_{1},\ldots,v_{n}\right)  ,$ $\left\vert
T\right\vert $ is the cardinality of $T$ and the vector $\frac{1}{\mathbf{p}%
}\mathbf{1}_{T}$ is defined by $\left(  \frac{1}{\mathbf{p}}\mathbf{1}%
_{T}\right)  _{i}=\frac{1}{p_{i}}$ if $i\in T$, $\left(  \frac{1}{\mathbf{p}%
}\mathbf{1}_{T}\right)  _{i}=0$ if $i\notin T$, for $i\in\left\{
1,2,\ldots,n\right\}  $.
\end{theorem}

This family is given, by example, for the simpler case $P\left(
\boldsymbol{\theta}\right)  =\prod_{i=1}^{n}(1+p_{i}\theta_{i})]-\beta
\mathbf{p}^{\left[  n\right]  }\boldsymbol{\theta}^{\left[  n\right]  },$ with
$0\leqslant\beta<1$, corresponding to the family gived by (Fang \textit{et
al.,} 2000), namely $C\left(  \mathbf{v}\right)  =\mathbf{v}^{\left[
n\right]  }[1-\beta\prod_{i=1}^{n}(1-v_{i}^{\frac{1}{\lambda}})]^{-\lambda}.$

\begin{theorem}
\label{copulagammaMM}Let $P$ an affine polynomial in the $n$ variables
$\theta_{i},$ $i=1,\ldots,n$, with $P\left(  \mathbf{0}\right)  =1,$ let
$\Lambda=\left(  \lambda,\lambda_{1},\ldots,\lambda_{n}\right)  $ with
$\lambda_{i}>\lambda>0,$ for $i\in\left\{  1,2,\ldots,n\right\}  ,$ then the
Laplace copula of the multivariate multi-factor gamma distribution
$\gamma_{\left(  P,\Lambda\right)  }$ associated to $\left(  P,\Lambda\right)
$ is%
\begin{equation}
C_{L_{\gamma_{\left(  P,\Lambda\right)  }}}\left(  \mathbf{v}\right)
=\mathbf{v}^{\left[  n\right]  }[1+\sum_{T\subset\left[  n\right]  ,\left\vert
T\right\vert >1}\left(  -1\right)  ^{\left\vert T\right\vert }P(-\frac
{1}{\mathbf{p}}\mathbf{1}_{T})\prod_{t\in T}(1-v_{t}^{\frac{1}{\lambda_{t}}%
})]^{-\lambda}. \label{LC_MMGD}%
\end{equation}

\end{theorem}

This copula family is a new generalization of Farlie-Gumbel-Morgenstern Copulas.

Bekrizadeh \textit{et al.,} 2012 propose a similar formula only for $\lambda$
being a negative integer.

We note that, if the conditions of Theorem \ref{mainresult} are checked,
formulas (\ref{LC_MGD}) and (\ref{LC_MMGD}) define copulas.

For these copulas, if we inject univariate gamma distributions $\gamma
_{\left(  p_{i},\lambda\right)  }$ in the first case or $\gamma_{\left(
p_{i},\lambda_{i}\right)  }$ in the second case, then we obtain other
multivariate gamma distributions and other multivariate multi-factor gamma
distributions in the sense that their marginal distributions are respectively
$\gamma_{\left(  p_{1},\lambda\right)  }$ and $\gamma_{\left(  p_{1}%
,\lambda_{i}\right)  }$. Their cdf and pdf are given by (\ref{Sklar1}) and
(\ref{Sklar2}) respectively, and have a link with \textit{multivariate gamma
distribution associated with} $\left(  P,\lambda\right)  $ and
\textit{multi-factor gamma distribution associated with} $\left(
P,\Lambda\right)  $ respectively.

\section{The bidimensional and tridimensional cases}

For the bidimensional case, Theorem (\ref{copulagammaM}) and Theorem
(\ref{copulagammaMM}) give the following corollary

\begin{corollary}
For the bivariate gamma distribution such that $L_{\gamma_{\left(
P,\lambda\right)  }}\left(  \boldsymbol{\theta}\right)  =\left(  P\left(
\boldsymbol{\theta}\right)  \right)  ^{-\lambda}$, we have
\[
C_{L_{\gamma_{\left(  P,\lambda\right)  }}}\left(  v_{1},v_{2}\right)
=v_{1}v_{2}[1-r_{12}(1-v_{1}^{\frac{1}{\lambda}})(1-v_{2}^{\frac{1}{\lambda}%
})]^{-\lambda},
\]
where, for $\mathbf{X}=\left(  X_{1},X_{2}\right)  $ with pd $\gamma_{\left(
P,\lambda\right)  },$ $r_{12}=-P(-p_{1}^{-1},-p_{2}^{-1})=1-p_{12}/\left(
p_{1}p_{2}\right)  $ is the linear correlation coefficient of the random
variables $X_{1},$ $X_{2}$ ; it checks $0\leqslant r_{12}\leqslant1$. It is
the copula of the BB10 family p. 154 in the Joe's book (Joe 1997). This result
is due to the following equalities : $E\left(  X_{i}\right)  =\lambda p_{i,}$
$Var\left(  X_{i}\right)  =\lambda p_{i}^{2},$ $Cov(X_{1},X_{2})=\lambda
\left(  -p_{12}+p_{1}p_{2}\right)  $ obtained by derivating $L_{\gamma
_{\left(  P,\lambda\right)  }}\left(  \boldsymbol{\theta}\right)  $ at
$\mathbf{0}.$\newline For the bivariate multi-factor gamma distribution such
that $L_{\gamma_{\left(  P,\Lambda\right)  }}\left(  \theta\right)  =\left(
P\left(  \theta\right)  \right)  ^{-\lambda}\prod_{i=1}^{2}\left(
1+p_{i}\theta_{i}\right)  ^{-\left(  \lambda_{i}-\lambda\right)  }$ , we have
a more general family
\begin{equation}
C_{L_{\gamma_{\left(  P,\Lambda\right)  }}}\left(  v_{1},v_{2}\right)
=v_{1}v_{2}[1-r_{12}(1-v_{1}^{\frac{1}{\lambda_{1}}})(1-v_{2}^{\frac
{1}{\lambda_{2}}})]^{-\lambda}, \label{bmfgLc}%
\end{equation}
where, for $\mathbf{X}=\left(  X_{1},X_{2}\right)  =\left(  Y_{1}%
,Y_{2}\right)  +\left(  Z_{1},Z_{2}\right)  =\mathbf{Y}+\mathbf{Z}$ as defined
in Proposition \ref{buildMFGD} with pd $\gamma_{\left(  P,\Lambda\right)  },$
and $0\leqslant r_{12}=-P(-p_{1}^{-1},-p_{2}^{-1})=1-p_{12}/\left(  p_{1}%
p_{2}\right)  \leqslant1,$ is the linear correlation coefficient of the random
variables $Y_{1},$ $Y_{2}$ in the bivariate gamma distribution such that
$L_{_{\mathbf{Y}}}\left(  \theta\right)  =\left[  P\left(  \theta\right)
\right]  ^{-\lambda}$.
\end{corollary}

This bivariate family is not given in the Joe's books (Joe, 1997, 2014).
Bekrizadeh \textit{et al.,} 2012 give a similar formula for $\lambda$ being a
negative integer.

We recall the formulas (Joe, 2014) for the computation of $\tau,$ the
Kendall's tau, and $\rho_{S},$ the Spearmann's rho%
\begin{align}
\tau &  =1-4%
{\displaystyle\iint_{\left[  0,1\right]  ^{2}}}
\frac{\partial C}{\partial u}\left(  u,v\right)  \frac{\partial C}{\partial
v}\left(  u,v\right)  \,\text{d}u\text{d}v,\label{tau formula}\\
\rho_{S}  &  =12%
{\displaystyle\iint_{\left[  0,1\right]  ^{2}}}
C\left(  u,v\right)  \,\text{d}u\text{d}v-3. \label{rho formula}%
\end{align}
Then, we can give the following result

\begin{proposition}
\label{Kendall's tau & Spearman's rho}The Kendall's tau and the Spearman's rho
of the Laplace copula of the bivariate multi-factor gamma distribution are
given by the following formulas%
\begin{align}
\tau &  =1-F_{2}^{3}\left(  2\lambda,1,1;2\lambda_{1}+1,2\lambda_{2}%
+1;r_{12}\right) \nonumber\\
&  +\frac{4\lambda}{\left(  2\lambda_{1}+1\right)  \left(  2\lambda
_{2}+1\right)  }r_{12}F_{2}^{3}\left(  2\lambda+1,1,2;2\lambda_{1}%
+2,2\lambda_{2}+2;r_{12}\right) \nonumber\\
&  -\frac{\lambda^{2}}{\left(  2\lambda_{1}+1\right)  \left(  2\lambda
_{2}+1\right)  \left(  \lambda_{1}+1\right)  \left(  \lambda_{2}+1\right)
}r_{12}^{2}F_{2}^{3}\left(  2\lambda+2,2,2;2\lambda_{1}+3,2\lambda
_{2}+3;r_{12}\right)  , \label{KtauHyper}%
\end{align}
and%
\begin{align}
\rho_{S}  &  =3\left[  F_{2}^{3}(1,1,\lambda;2\lambda_{1}+1,2\lambda
_{2}+1;r_{12})-1\right] \label{SrhoHyper1}\\
&  =\frac{3\lambda}{\left(  2\lambda_{1}+1\right)  \left(  2\lambda
_{2}+1\right)  }r_{12}F_{2}^{3}(1,2,\lambda+1;2\lambda_{1}+2,2\lambda
_{2}+2;r_{12}). \label{SrhoHyper2}%
\end{align}

\end{proposition}

For the tridimensional case, Theorem (\ref{copulagammaM}) and Theorem
(\ref{copulagammaMM}) give the following corollary

\begin{corollary}
For the trivariate gamma distribution $\gamma_{\left(  P,\lambda\right)  }$
such that $L_{\gamma_{\left(  P,\lambda\right)  }}\left(  \boldsymbol{\theta
}\right)  =\left(  P\left(  \boldsymbol{\theta}\right)  \right)  ^{-\lambda}$,
we have
\begin{align*}
C_{L_{\gamma_{\left(  P,\lambda\right)  }}}\left(  v_{1},v_{2},v_{3}\right)
&  =v_{1}v_{2}v_{3}[1-r_{12}(1-v_{1}^{\frac{1}{\lambda}})(1-v_{2}^{\frac
{1}{\lambda}})-r_{13}(1-v_{1}^{\frac{1}{\lambda}})(1-v_{3}^{\frac{1}{\lambda}%
})\\
&  -r_{23}(1-v_{2}^{\frac{1}{\lambda}})(1-v_{3}^{\frac{1}{\lambda}}%
)+2r_{123}(1-v_{1}^{\frac{1}{\lambda}})(1-v_{2}^{\frac{1}{\lambda}}%
)(1-v_{3}^{\frac{1}{\lambda}})]^{-\lambda},
\end{align*}
\newline where, for $\mathbf{X}=\left(  X_{1},X_{2},X_{3}\right)  $ with pd
$\gamma_{\left(  P,\lambda\right)  },$ we denote by $r_{ij}=-P(-\frac{1}%
{p}\mathbf{1}_{\left\{  i,j\right\}  })=1-p_{ij}/(p_{i}p_{j}),$ $1\leqslant
i\neq j\leqslant3,$ the linear correlation coefficient of the random variables
$X_{i},X_{j},$ they check $0\leqslant r_{ij}\leqslant1,$ and we denote by
$r_{123}$ the number defined by
\begin{align}
r_{123}  &  =\mathbb{E}\{\prod_{i=1}^{3}\left[  X_{i}-\mathbb{E}\left(
X_{i}\right)  \right]  \}/\prod_{i=1}^{3}(\mathbb{E}\{\left[  X_{i}%
-\mathbb{E}\left(  X_{i}\right)  \right]  ^{3}\})^{1/3}\label{r123}\\
&  =-\tfrac{1}{2}P(-p_{1}^{-1},-p_{2}^{-1},-p_{3}^{-1})\nonumber
\end{align}
(to compare with $r_{12}=\mathbb{E}\{\prod_{i=1}^{2}\left[  X_{i}%
-\mathbb{E}\left(  X_{i}\right)  \right]  \}/\prod_{i=1}^{2}(\mathbb{E}%
\{\left[  X_{i}-\mathbb{E}\left(  X_{i}\right)  \right]  ^{2}\})^{1/2}%
$).\newline This last result is due to the following equalities :
$\mathbb{E}\left(  X_{i}\right)  =\lambda p_{i,}$ $\mathbb{E}\{\left[
X_{i}-\mathbb{E}\left(  X_{i}\right)  \right]  ^{3}\}=2\lambda p_{i}^{3},$
\newline$\mathbb{E}\{\prod_{i=1}^{3}\left[  X_{i}-\mathbb{E}\left(
X_{i}\right)  \right]  \}=-\lambda p_{1}p_{2}p_{3}P(-p_{1}^{-1},-p_{2}%
^{-1},-p_{3}^{-1})$ obtained by derivating $L_{\gamma_{\left(  P,\lambda
\right)  }}\left(  \boldsymbol{\theta}\right)  $ at $\mathbf{0}.$\newline For
the trivariate multi-factor gamma distribution such that $L_{\gamma_{\left(
P,\Lambda\right)  }}\left(  \boldsymbol{\theta}\right)  =\left[  P\left(
\boldsymbol{\theta}\right)  \right]  ^{-\lambda}\prod_{i=1}^{3}\left(
1+p_{i}\theta_{i}\right)  ^{-\left(  \lambda_{i}-\lambda\right)  }$
($\lambda_{i}\geqslant\lambda,$ $i=1,2,3$ ), we have%
\begin{align*}
C_{L_{\gamma_{\left(  P,\Lambda\right)  }}}\left(  v_{1},v_{2},v_{3}\right)
&  =v_{1}v_{2}v_{3}[1-r_{12}(1-v_{1}^{\frac{1}{\lambda_{1}}})(1-v_{2}%
^{\frac{1}{\lambda_{2}}})-r_{13}(1-v_{1}^{\frac{1}{\lambda_{1}}}%
)(1-v_{3}^{\frac{1}{\lambda_{3}}})\\
&  -r_{23}(1-v_{2}^{\frac{1}{\lambda_{2}}})(1-v_{3}^{\frac{1}{\lambda_{3}}%
})+2r_{123}(1-v_{1}^{\frac{1}{\lambda_{1}}})(1-v_{2}^{\frac{1}{\lambda_{2}}%
})(1-v_{3}^{\frac{1}{\lambda_{3}}})]^{-\lambda},
\end{align*}
where, for $\mathbf{X}=\left(  X_{1},X_{2},X_{3}\right)  =\left(  Y_{1}%
,Y_{2},Y_{3}\right)  +\left(  Z_{1},Z_{2},Z_{3}\right)  =\mathbf{Y}%
+\mathbf{Z}$ as defined in Proposition \ref{buildMFGD} with pd $\gamma
_{\left(  P,\Lambda\right)  },$ we denote by $r_{ij}=-P(-\frac{1}{\mathbf{p}%
}\mathbf{1}_{\left\{  i,j\right\}  })=1-\frac{p_{ij}}{p_{i}p_{j}},$
$1\leqslant i\neq j\leqslant3,$ the linear correlation coefficient of the
random variables $Y_{i},$ $Y_{j}$ in the trivariate gamma distribution such
that $L_{_{Y}}\left(  \boldsymbol{\theta}\right)  =\left(  P\left(
\boldsymbol{\theta}\right)  \right)  ^{-\lambda},$ they check $0\leqslant
r_{ij}\leqslant1,$ and $r_{123}$ is defined by equality (\ref{r123}).
\end{corollary}

\section{Appendix}

\subsection{Proof of Theorem \ref{Proposition_multifactorbig2}}

Let $\mathbf{X}=(X_{1},X_{2})$ and $\mathbf{Y}=(Y_{1},Y_{2})$ be independent
random variables. The pdf of $\mathbf{Z}=(X_{1}+Y_{1},X_{2}+Y_{2})$ where
$X_{i}$ has pd $\gamma\left(  p_{i},\lambda_{i}-\lambda\right)  ,$ $i=1,2$ and
$\mathbf{Y}$ has pd $\gamma\left(  P,\lambda\right)  $, is obtained by
convolution. By changing variables $1-v_{i}/z_{i}=u_{i}$, $i=1,2,$ and with
the notation $c=b_{12}(\widetilde{P})=-p_{12}^{-1}+p_{1}p_{2}p_{12}^{-2},$
$d_{i}=cp_{12}p_{i}^{-1},$ $i=1,2$ we obtain the pdf $f_{\mathbf{Z}}$ of
$\mathbf{Z}$ by $f_{\boldsymbol{Z}}\left(  z_{1},z_{2}\right)  =0$ if $\left(
z_{1},z_{2}\right)  \notin\left(  0,\infty\right)  ^{2}$ and for $\left(
z_{1},z_{2}\right)  \in\left(  0,\infty\right)  ^{2}$%
\begin{align}
f_{\mathbf{Z}}\left(  z_{1},z_{2}\right)   &  =\int_{0}^{z_{1}}\int_{0}%
^{z_{2}}\frac{1}{\left[  \Gamma\left(  \lambda\right)  \right]  ^{2}%
p_{12}^{\lambda}}\mathbf{e}^{-\left(  \frac{p_{2}}{p_{12}}v_{1}+\frac{p_{1}%
}{p_{12}}v_{2}\right)  }\left(  v_{1}v_{2}\right)  ^{\lambda-1}F_{1}\left(
\lambda;cv_{1}v_{2}\right) \nonumber\\
&  \times\frac{\mathbf{e}^{-\frac{z_{1}-v_{1}}{p_{1}}}}{\Gamma\left(
\lambda_{1}-\lambda\right)  }\left(  \frac{z_{1}-v_{1}}{p_{1}}\right)
^{\lambda_{1}-\lambda}\frac{\mathbf{e}^{-\frac{z_{2}-v_{2}}{p_{2}}}}%
{\Gamma\left(  \lambda_{2}-\lambda\right)  }\left(  \frac{z_{2}-v_{2}}{p_{2}%
}\right)  ^{\lambda_{2}-\lambda}\frac{\text{d}v_{1}}{z_{1}-v_{1}}%
\frac{\text{d}v_{2}}{z_{2}-v_{2}}\nonumber\\
&  =\frac{z_{1}^{\lambda_{1}-1}z_{2}^{\lambda_{2}-1}\mathbf{e}^{-\left(
\frac{p_{2}}{p_{12}}z_{1}+\frac{p_{1}}{p_{12}}z_{2}\right)  }}{\left(
\Gamma\left(  \lambda\right)  \right)  \Gamma\left(  \lambda_{1}%
-\lambda\right)  \Gamma\left(  \lambda_{2}-\lambda\right)  p_{12}^{\lambda
}p_{1}^{\lambda_{1}-\lambda}p_{2}^{\lambda_{2}-\lambda}}\sum_{k\geqslant
0}\frac{\left(  cz_{1}z_{2}\right)  ^{k}}{\Gamma\left(  \lambda+k\right)
k!}\nonumber\\
&  \times\int_{0}^{1}\mathbf{e}^{d_{1}z_{1}u_{1}}u_{1}^{\lambda_{1}-\lambda
-1}\left(  1-u_{1}\right)  ^{\lambda+k-1}\,\text{d}u_{1}\int_{0}^{1}%
\mathbf{e}^{d_{2}z_{2}u_{2}}u_{2}^{\lambda_{2}-\lambda-1}\left(
1-u_{2}\right)  ^{\lambda+k-1}\,\text{d}u_{2}. \label{convolution_gamma_b}%
\end{align}
As we have, with the notation $B\left(  \alpha,\beta\right)  =\int_{0}%
^{1}u^{\alpha-1}\left(  1-u\right)  ^{\beta-1}\,$d$u=\Gamma\left(
\alpha\right)  \Gamma\left(  \beta\right)  /\Gamma\left(  \alpha+\beta\right)
$, $\alpha,\beta>0$ for the Euler's Beta function,%
\[
\int_{0}^{1}\mathbf{e}^{\delta u}u^{\alpha-1}\left(  1-u\right)  ^{\beta
-1}\,\text{d}u=B\left(  \alpha,\beta\right)  \sum_{n=0}^{\infty}\frac{\left(
\alpha\right)  _{n}}{\left(  \alpha+\beta\right)  _{n}}\frac{\delta^{n}}{n!},
\]
we obtain%
\begin{align*}
f_{\mathbf{Z}}\left(  \boldsymbol{z}\right)   &  =\frac{z_{1}^{\lambda_{1}%
-1}z_{2}^{\lambda_{2}-1}\mathbf{e}^{-\left(  \frac{p_{2}}{p_{12}}z_{1}%
+\frac{p_{1}}{p_{12}}z_{2}\right)  }}{\Gamma\left(  \lambda\right)
\Gamma\left(  \lambda_{1}-\lambda\right)  \Gamma\left(  \lambda_{2}%
-\lambda\right)  p_{12}^{\lambda}p_{1}^{\lambda_{1}-\lambda}p_{2}^{\lambda
_{2}-\lambda}}\sum_{k\geqslant0}\frac{\left(  cz_{1}z_{2}\right)  ^{k}}%
{\Gamma\left(  \lambda+k\right)  k!}\frac{\Gamma\left(  \lambda_{1}%
-\lambda\right)  \Gamma\left(  \lambda+k\right)  }{\Gamma\left(  \lambda
_{1}+k\right)  }\\
&  \times\sum_{n=0}^{\infty}\frac{\left(  \lambda_{1}-\lambda\right)  _{n}%
}{\left(  \lambda_{1}+k\right)  _{n}}\frac{\left(  d_{1}z_{1}\right)  ^{n}%
}{n!}\frac{\Gamma\left(  \lambda_{2}-\lambda\right)  \Gamma\left(
\lambda+k\right)  }{\Gamma\left(  \lambda_{2}+k\right)  }\sum_{m=0}^{\infty
}\frac{\left(  \lambda_{2}-\lambda\right)  _{m}}{\left(  \lambda_{2}+k\right)
_{m}}\frac{\left(  d_{2}z_{2}\right)  ^{m}}{m!}\\
&  =\frac{z_{1}^{\lambda_{1}-1}z_{2}^{\lambda_{2}-1}\mathbf{e}^{-\left(
\frac{p_{2}}{p_{12}}z_{1}+\frac{p_{1}}{p_{12}}z_{2}\right)  }}{\Gamma\left(
\lambda_{1}\right)  \Gamma\left(  \lambda_{2}\right)  p_{12}^{\lambda}%
p_{1}^{\lambda_{1}-\lambda}p_{2}^{\lambda_{2}-\lambda}}\sum_{k=0}^{\infty}%
\sum_{n=0}^{\infty}\sum_{m=0}^{\infty}\frac{\left(  \lambda_{1}-\lambda
\right)  _{n}\left(  \lambda_{2}-\lambda\right)  _{m}\left(  \lambda\right)
_{k}}{\left(  \lambda_{1}\right)  _{k+n}\left(  \lambda_{2}\right)  _{k+m}%
}\frac{\left(  d_{1}z_{1}\right)  ^{n}}{n!}\frac{\left(  d_{2}z_{2}\right)
^{m}}{m!}\frac{\left(  cz_{1}z_{2}\right)  ^{k}}{k!}.
\end{align*}
Hence, we have proved the formula (\ref{mulbigamma2}).

\subsection{Proof of Theorem \ref{gamma3}}

We start from the Laplace transform of the trivariate gamma distribution
$\gamma_{\left(  \lambda,P\right)  }$ associated to $\left(  \lambda,P\right)
$. First, we write%

\begin{align*}
L_{\gamma_{\left(  \lambda,P\right)  }}\left(  \boldsymbol{\theta}\right)   &
=\left(  1+p_{1}\theta_{1}+p_{2}\theta_{2}+p_{3}\theta_{3}+p_{12}\theta
_{1}\theta_{2}+p_{13}\theta_{1}\theta_{3}+p_{23}\theta_{2}\theta_{3}%
+p_{123}\theta_{1}\theta_{2}\theta_{3}\right)  ^{-\lambda}\\
&  =\left(  1+p_{3}\theta_{3}\right)  ^{-\lambda}\left(  1+\frac{p_{1}%
+p_{13}\theta_{3}}{1+p_{3}\theta_{3}}\theta_{1}+\frac{p_{2}+p_{23}\theta_{3}%
}{1+p_{3}\theta_{3}}\theta_{2}+\frac{p_{12}+p_{123}\theta_{3}}{1+p_{3}%
\theta_{3}}\theta_{1}\theta_{2}\right)  ^{-\lambda}.
\end{align*}
Let $Q\left(  \theta_{1},\theta_{2}\right)  =1+q_{1}\theta_{1}+q_{2}\theta
_{2}+q_{12}\theta_{1}\theta_{2}$ where $q_{1}=\left(  p_{1}+p_{13}\theta
_{3}\right)  /\left(  1+p_{3}\theta_{3}\right)  ,$ $q_{2}=\left(  p_{2}%
+p_{23}\theta_{3}\right)  /\left(  1+p_{3}\theta_{3}\right)  ,$ and
$q_{12}=\left(  p_{12}+p_{123}\theta_{3}\right)  /\left(  1+p_{3}\theta
_{3}\right)  $. With these notations, we have
\[
L_{\gamma_{\left(  \lambda,P\right)  }}\left(  \boldsymbol{\theta}\right)
=\left(  1+p_{3}\theta_{3}\right)  ^{-\lambda}\left(  1+q_{1}\theta_{1}%
+q_{2}\theta_{2}+q_{12}\theta_{1}\theta_{2}\right)  ^{-\lambda}.
\]
\newline Let $\widetilde{p}_{T}=-p_{\overline{T}}/p_{123}$ ,$\widetilde{b}%
_{ij}=\widetilde{p}_{ij}+\widetilde{p}_{i}\widetilde{p}_{j},$ and
$\widetilde{b}_{123}=\widetilde{p}_{123}+\widetilde{p}_{1}\widetilde{p}%
_{23}+\widetilde{p}_{2}\widetilde{p}_{13}+\widetilde{p}_{3}\widetilde{p}%
_{12}+2\widetilde{p}_{1}\widetilde{p}_{2}\widetilde{p}_{3},$ then we have%
\begin{align*}
b_{1,2}(\widetilde{Q})  &  =\frac{-q_{12}+q_{1}q_{2}}{q_{12}^{2}}\\
&  =\widetilde{b}_{12}+\frac{\widetilde{b}_{123}}{\left(  \theta
_{3}-\widetilde{p}_{3}\right)  }+\frac{\widetilde{b}_{13}\widetilde{b}_{23}%
}{\left(  \theta_{3}-\widetilde{p}_{3}\right)  ^{2}}.
\end{align*}

Let $\gamma_{\left(  \lambda,Q\right)  }$ be the bivariate gamma distribution
associated to $\left(  \lambda,Q\right)  $, its Laplace transform
is\newline$\left(  1+q_{1}\theta_{1}+q_{2}\theta_{2}+q_{12}\theta_{2}%
\theta_{3}\right)  ^{-\lambda}$. We have%
\[
\gamma_{\left(  \lambda,Q\right)  }\left(  dx\right)  =\frac{q_{12}^{-\lambda
}}{\left[  \Gamma\left(  \lambda\right)  \right]  ^{2}}e^{-\frac{q_{2}}%
{q_{12}}x_{1}-\frac{q_{1}}{q_{12}}x_{2}}\left(  x_{1}x_{2}\right)
^{\lambda-1}F_{1}[\lambda,b_{1,2}(\widetilde{Q})x_{1}x_{2}%
]\mathbbm{1}_{\left(  0,\infty\right)  ^{2}}\left(  x_{1},x_{2}\right)
\,\text{d}x_{1}\text{d}x_{2}.
\]
Second, we are looking for the unidimensional distribution for which its
Laplace transform with respect to the variable $\theta_{3}$ is equal to%
\begin{align}
&  \left(  1+p_{3}\theta_{3}\right)  ^{-\lambda}\frac{q_{12}^{-\lambda}%
}{\left(  \Gamma\left(  \lambda\right)  \right)  ^{2}}e^{-\frac{q_{2}}{q_{12}%
}x_{1}-\frac{q_{1}}{q_{12}}x_{2}}\left(  x_{1}x_{2}\right)  ^{\lambda-1}%
F_{1}[\lambda,b_{1,2}(\widetilde{Q})x_{1}x_{2}]\mathbbm{1}_{\left(
0,\infty\right)  ^{2}}\left(  x_{1},x_{2}\right)  \,\text{d}x_{1}\text{d}%
x_{2}\nonumber\\
&  =\frac{p_{123}^{-\lambda}}{\left[  \Gamma\left(  \lambda\right)  \right]
^{2}}\exp\left(  \widetilde{p}_{1}x_{1}+\widetilde{p}_{2}x_{2}\right)  \left(
x_{1}x_{2}\right)  ^{\lambda-1}\left(  \theta_{3}-\widetilde{p}_{3}\right)
^{-\lambda}\exp(\frac{\widetilde{b}_{13}x_{1}+\widetilde{b}_{23}x_{2}}%
{\theta_{3}-\widetilde{p}_{3}})\nonumber\\
&  \times F_{1}\{\lambda,[\widetilde{b}_{12}+\frac{\widetilde{b}_{123}}%
{\theta_{3}-\widetilde{p}_{3}}+\frac{\widetilde{b}_{13}\widetilde{b}_{23}%
}{\left(  \theta_{3}-\widetilde{p}_{3}\right)  ^{2}}]x_{1}x_{2}%
\}\mathbbm{1}_{\left(  0,\infty\right)  ^{2}}\left(  x_{1},x_{2}\right)
\,\text{d}x_{1}\text{d}x_{2}. \label{MGD3STEP1}%
\end{align}

From equality (\ref{MGD3STEP1}) and (Hladik, 1969), we obtain the
unidimensional distribution for which its Laplace transform is%
\begin{align}
&  \frac{p_{123}^{-\lambda}}{\left[  \Gamma\left(  \lambda\right)  \right]
^{2}}\exp\left(  \widetilde{p}_{1}x_{1}+\widetilde{p}_{2}x_{2}+\widetilde{p}%
_{3}x_{3}\right)  \left(  x_{1}x_{2}\right)  ^{\lambda-1}\theta_{3}^{-\lambda
}\exp(\frac{\widetilde{b}_{13}x_{1}+\widetilde{b}_{23}x_{2}}{\theta_{3}%
})\nonumber\\
&  \times F_{1}[\lambda,(\frac{\widetilde{b}_{13}\widetilde{b}_{23}}%
{\theta_{3}^{2}}+\frac{\widetilde{b}_{123}}{\theta_{3}}+\widetilde{b}%
_{12})x_{1}x_{2}]\mathbbm{1}_{\left(  0,\infty\right)  ^{2}}\left(
x_{1},x_{2}\right)  \,\text{d}x_{1}\text{d}x_{2}. \label{MGDSTEP1B}%
\end{align}

Third, we are looking for the unidimensional distribution for which its
Laplace transform with respect to the variable $\theta_{3}$ is equal to%
\begin{align}
&  \theta_{3}^{-\lambda}\exp(\frac{\widetilde{b}_{13}x_{1}+\widetilde{b}%
_{23}x_{2}}{\theta_{3}})F_{1}[\lambda,(\frac{\widetilde{b}_{13}\widetilde{b}%
_{23}}{\theta_{3}^{2}}+\frac{\widetilde{b}_{123}}{\theta_{3}}+\widetilde{b}%
_{12})x_{1}x_{2}]\nonumber\\
&  =\sum_{k=0}^{\infty}\frac{1}{\left(  \lambda\right)  _{k}}\left(
x_{1}x_{2}\right)  ^{k}\sum_{\ell+m+n=k}\frac{1}{\ell!m!n!}\widetilde{b}%
_{13}^{\ell}\widetilde{b}_{23}^{\ell}\widetilde{b}_{123}^{m}\widetilde{b}%
_{12}^{n}\theta_{3}^{-\left(  \lambda+2\ell+m\right)  }\exp(\frac
{\widetilde{b}_{13}x_{1}+\widetilde{b}_{23}x_{2}}{\theta_{3}})
\label{MGD3STEP2}%
\end{align}
From (Hadlik,1986), we have the following equality%
\begin{equation}
L_{\left[  \Gamma\left(  \lambda\right)  \right]  ^{-1}t^{\lambda-1}%
F_{1}\left(  \lambda,at\right)  \mathbbm{1}_{\left(  0,\infty\right)  }\left(
t\right)  \,\text{d}t}\left(  s\right)  =s^{-\lambda}\exp\left(  \frac{a}%
{s}\right)  \label{hladik}%
\end{equation}
We utilize the equality (\ref{hladik}) in (\ref{MGD3STEP2}) and we obtain the
following unidimensional distribution%

\begin{align}
&  \sum_{k=0}^{\infty}\frac{1}{\left(  \lambda\right)  _{k}}\left(  x_{1}%
x_{2}\right)  ^{k}\sum_{\ell+m+n=k}\frac{1}{\ell!m!n!}\widetilde{b}_{13}%
^{\ell}\widetilde{b}_{23}^{\ell}\widetilde{b}_{123}^{m}\widetilde{b}_{12}%
^{n}\left[  \Gamma\left(  \lambda+2\ell+m\right)  \right]  ^{-1}x_{3}%
^{\lambda+2\ell+m-1}\nonumber\\
&  \times F_{1}\left[  \lambda+2\ell+m,(\widetilde{b}_{13}x_{1}+\widetilde{b}%
_{23}x_{2})x_{3}\right]  \mathbbm{1}_{\left(  0,\infty\right)  }\left(
x_{3}\right)  \,\text{d}x_{3}\nonumber\\
&  =\left[  \Gamma\left(  \lambda\right)  \right]  ^{-1}x_{3}^{\lambda
-1}\mathbbm{1}_{\left(  0,\infty\right)  }\left(  x_{3}\right) \nonumber\\
&  \times\sum_{\left(  \ell,m,n\right)  \in\mathbb{N}^{3}}\frac{1}{\left(
\lambda\right)  _{\left(  \ell+m+n\right)  }\left(  \lambda\right)  _{\left(
2\ell+m\right)  }\ell!m!n!}(\widetilde{b}_{13}x_{1}x_{3}\widetilde{b}%
_{23}x_{2}x_{3})^{\ell}(\widetilde{b}_{123}x_{1}x_{2}x_{3})^{m}(\widetilde{b}%
_{12}x_{1}x_{2})^{n}\nonumber\\
&  \times F_{1}(\lambda+2\ell+m,\widetilde{b}_{13}x_{1}x_{3}+\widetilde{b}%
_{23}x_{2}x_{3})\,\text{d}x_{3}\nonumber\\
&  =\left[  \Gamma\left(  \lambda\right)  \right]  ^{-1}x_{3}^{\lambda
-1}\mathbbm{1}_{\left(  0,\infty\right)  }\left(  x_{3}\right) \nonumber\\
&  \times\sum_{\left(  \ell,m,n,p\right)  \in\mathbb{N}^{4}}\frac{1}{\left(
\lambda\right)  _{\left(  \ell+m+n\right)  }\left(  \lambda\right)  _{\left(
2\ell+m+p\right)  }\ell!m!n!p!}\nonumber\\
&  \times(\widetilde{b}_{13}x_{1}x_{3}\widetilde{b}_{23}x_{2}x_{3})^{\ell
}(\widetilde{b}_{123}x_{1}x_{2}x_{3})^{m}(\widetilde{b}_{12}x_{1}x_{2}%
)^{n}(\widetilde{b}_{13}x_{1}x_{3}+\widetilde{b}_{23}x_{2}x_{3})^{p}%
\,\text{d}x_{3}\nonumber\\
&  =\left[  \Gamma\left(  \lambda\right)  \right]  ^{-1}x_{3}^{\lambda
-1}F_{II}(\lambda,\lambda,\widetilde{b}_{13}x_{1}x_{3}\widetilde{b}_{23}%
x_{2}x_{3},\widetilde{b}_{123}x_{1}x_{2}x_{3},\widetilde{b}_{12}x_{1}%
x_{2},\widetilde{b}_{13}x_{1}x_{3}+\widetilde{b}_{23}x_{2}x_{3}%
)\mathbbm{1}_{\left(  0,\infty\right)  }\left(  x_{3}\right)  \,\text{d}x_{3}
\label{MGD3STEP3}%
\end{align}
We carry this last equality (\ref{MGD3STEP3}) in (\ref{MGDSTEP1B}) and we
obtain the pd (\ref{MGD3}).

\subsection{Proof of Theorem \ref{copulagammaM}}

We utilize the following notations, for $i\in\left[  n\right]  ,$
$L_{\gamma_{\left(  p_{i},\lambda\right)  }}\left(  \theta_{i}\right)
=\left(  1+p_{i}\theta_{i}\right)  ^{-\lambda}=v_{i},0<\lambda_{i}<\lambda$
and $1+p_{i}\theta_{i}=u_{i}$. In particular $\theta_{i}=\left(
u_{i}-1\right)  /p_{i}=\left(  v_{i}^{-1/\lambda}-1\right)  /p_{i}$ and
$u_{i}=v_{i}^{-1/\lambda}.$ We remark that $u_{i}=0\Leftrightarrow\theta
_{i}=-1/p_{i}.$

We start from the Laplace transform of the \textit{multivariate gamma
distribution }associated to $\left(  P,\lambda\right)  $. From the definition
of $\gamma_{\left(  P,\lambda\right)  }$, we obtain%
\[
\lbrack L_{\gamma_{\left(  P,\lambda\right)  }}\left(  \boldsymbol{\theta
}\right)  ]^{-\frac{1}{\lambda}}=\sum_{T\subset\left[  n\right]  }%
p_{T}\boldsymbol{\theta}^{T}.
\]
Hence, we have%
\begin{align*}
\lbrack L_{\gamma_{\left(  P,\lambda\right)  }}\left(  \boldsymbol{\theta
}\right)  ]^{-\frac{1}{\lambda}}  &  =\sum_{T\subset\left[  n\right]  }%
p_{T}\left(  \mathbf{u}-\mathbf{1}_{\left[  n\right]  }\right)  ^{T}\left(
\frac{1}{\mathbf{p}}\right)  ^{T}\\
&  =\sum_{T\subset\left[  n\right]  }q_{T}\mathbf{u}^{T}=Q\left(
\mathbf{u}\right)
\end{align*}
where $\mathbf{u}=\left(  u_{1},\ldots,u_{n}\right)  ,$ and
\[
Q\left(  \mathbf{u}\right)  =\sum_{T\subset\left[  n\right]  }\alpha
_{T}\left(  \mathbf{u}-\mathbf{1}_{\left[  n\right]  }\right)  ^{T}%
\mathbf{u}^{\overline{T}},
\]
indeed, the polynomials $\left(  \mathbf{u}-\mathbf{1}_{\left[  n\right]
}\right)  ^{T}\mathbf{u}^{\overline{T}},$ $T\subset\left[  n\right]  $ form a
basis of the vector space of affine polynomials of degree less or equal to
$n.$ Because $u_{i}=0\Leftrightarrow\theta_{i}=-\frac{1}{p_{i}},$ for
$i\in\left[  n\right]  ,$ we have%
\[
Q\left(  0\right)  =P(-\frac{1}{\mathbf{p}}\mathbf{1}_{\left[  n\right]
})=\left(  -1\right)  ^{n}\alpha_{\left[  n\right]  },
\]
consequently
\[
\alpha_{\left[  n\right]  }=\left(  -1\right)  ^{n}P(-\frac{1}{\mathbf{p}%
}\mathbf{1}_{\left[  n\right]  }).
\]
Since $P$ is an affine polynomial, we have
\[
\alpha_{T}=\left(  -1\right)  ^{\left\vert T\right\vert }P(-\frac
{1}{\mathbf{p}}\mathbf{1}_{T}),
\]
in particular $\alpha_{T}=0,$ if $\left\vert T\right\vert =1$ and
$\alpha_{\varnothing}=1.$ Finally we obtain%
\begin{align}
P\left(  \boldsymbol{\theta}\right)   &  =\sum_{T\subset\left[  n\right]
}\left(  -1\right)  ^{\left\vert T\right\vert }P(-\frac{1}{\mathbf{p}%
}\mathbf{1}_{T})\left(  \mathbf{u}-\mathbf{1}_{\left[  n\right]  }\right)
^{T}\mathbf{u}^{\overline{T}}\label{d_p_affine}\\
&  =\mathbf{u}^{\left[  n\right]  }+\sum_{T\subset\left[  n\right]
,\left\vert T\right\vert >1}\left(  -1\right)  ^{\left\vert T\right\vert
}P(-\frac{1}{\mathbf{p}}\mathbf{1}_{T})\left(  \mathbf{u}-\mathbf{1}_{\left[
n\right]  }\right)  ^{T}\mathbf{u}^{\overline{T}}\nonumber
\end{align}
We deduce from equality (\ref{d_p_affine})%
\begin{align*}
L_{\mathbf{\gamma}_{\left(  P,\lambda\right)  }}\left(  \boldsymbol{\theta
}\right)   &  =[\sum_{T\subset\left[  n\right]  }\left(  -1\right)
^{\left\vert T\right\vert }P(-\frac{1}{\mathbf{p}}\mathbf{1}_{T})\left(
\mathbf{u}-\mathbf{1}_{\left[  n\right]  }\right)  ^{T}\mathbf{u}%
^{\overline{T}}]^{-\lambda}\\
&  =\left(  \mathbf{u}^{\left[  n\right]  }\right)  ^{-\lambda}[\sum
_{T\subset\left[  n\right]  }\left(  -1\right)  ^{\left\vert T\right\vert
}P(-\frac{1}{\mathbf{p}}\mathbf{1}_{T})(\mathbf{1}_{\left[  n\right]  }%
-\frac{1}{\mathbf{u}})^{T}]^{-\lambda}\\
&  =\mathbf{v}^{\left[  n\right]  }[1+\sum_{T\subset\left[  n\right]
,\left\vert T\right\vert >1}\left(  -1\right)  ^{\left\vert T\right\vert
}P(-\frac{1}{\mathbf{p}}\mathbf{1}_{T})\prod_{t\in T}(1-v_{t}^{\frac
{1}{\lambda}})]^{-\lambda}\\
&  =C_{L_{\gamma_{\left(  P,\lambda\right)  }}}\left(  \mathbf{v}\right)  .
\end{align*}

\subsection{Proof of Theorem \ref{copulagammaMM}}

We start from the Laplace transform of the \textit{multivariate multi-factor
gamma distribution }associated to $\left(  P,\Lambda\right)  ,$ $\Lambda
=\left(  \lambda,\lambda_{1},\ldots\lambda_{n}\right)  $, where $P$ is an
affine polynomial with respect to $\boldsymbol{\theta}=\left(  \theta
_{1},\ldots,\theta_{n}\right)  $ and $\lambda_{i}\geqslant\lambda,$ for
$i\in\left[  n\right]  .$ We use the following notations $u_{i}=\left(
1+p_{i}\theta_{i}\right)  ^{-\lambda}$ and $v_{i}=\left(  1+p_{i}\theta
_{i}\right)  ^{-\lambda_{i}}$ so that $\left(  1+p_{i}\theta_{i}\right)
^{-\left(  \lambda_{i}-\lambda\right)  }=v_{i}/u_{i}$ and $\left(
1+p_{i}\theta_{i}\right)  ^{-1}=u_{i}^{1/\lambda}=v_{i}^{1/\lambda_{i}}.$
According to equality \ref{LC_MGD}, we have%
\begin{align*}
L_{\mathbf{\gamma}_{\left(  P,\Lambda\right)  }}\left(  \boldsymbol{\theta
}\right)   &  =\left[  P\left(  \boldsymbol{\theta}\right)  \right]
^{-\lambda}\prod_{i\in\left[  n\right]  }\left(  1+p_{i}\theta_{i}\right)
^{-\left(  \lambda_{i}-\lambda\right)  }\\
&  =\mathbf{u}^{\left[  n\right]  }[1+\sum_{T\subset\left[  n\right]
,\left\vert T\right\vert >1}\left(  -1\right)  ^{\left\vert T\right\vert
}P(-\frac{1}{\mathbf{p}}\mathbf{1}_{T})\prod_{t\in T}(1-u_{t}^{\frac
{1}{\lambda}})]^{-\lambda}\prod_{i\in\left[  n\right]  }\frac{v_{i}}{u_{i}}\\
&  =\mathbf{v}^{\left[  n\right]  }[1+\sum_{T\subset\left[  n\right]
,\left\vert T\right\vert >1}\left(  -1\right)  ^{\left\vert T\right\vert
}P(-\frac{1}{\mathbf{p}}\mathbf{1}_{T})\prod_{t\in T}(1-u_{t}^{\frac
{1}{\lambda}})]^{-\lambda}\\
&  =\mathbf{v}^{\left[  n\right]  }[1+\sum_{T\subset\left[  n\right]
,\left\vert T\right\vert >1}\left(  -1\right)  ^{\left\vert T\right\vert
}P(-\frac{1}{\mathbf{p}}\mathbf{1}_{T})\prod_{t\in T}(1-v_{t}^{\frac
{1}{\lambda_{t}}})]^{-\lambda}\\
&  =C_{L_{\gamma_{\left(  P,\Lambda\right)  }}}\left(  \mathbf{v}\right)  .
\end{align*}

\subsection{Proof of Proposition \ref{Kendall's tau & Spearman's rho}}

\subsubsection{Kendall's tau}

By injecting the given copula in formula (\ref{bmfgLc}) in formula
(\ref{tau formula}), we obtain%

\begin{align*}
\frac{1-\tau}{4}  &  =%
{\displaystyle\iint_{\left[  0,1\right]  ^{2}}}
u_{1}u_{2}[1-r_{12}(1-u_{1}^{\frac{1}{\lambda_{1}}})(1-u_{2}^{\frac{1}%
{\lambda_{2}}})]^{-2\lambda-2}\\
&  \{1-r_{12}[1-(1-\frac{\lambda}{\lambda_{1}})u_{1}^{\frac{1}{\lambda_{1}}%
}](1-u_{2}^{\frac{1}{\lambda_{2}}})\}\{1-r_{12}[1-(1-\frac{\lambda}%
{\lambda_{2}})u_{2}^{\frac{1}{\lambda_{2}}}](1-u_{1}^{\frac{1}{\lambda_{1}}%
})\}\,\text{d}u_{1}\text{d}u_{2}.
\end{align*}
\newline By changing the variables $t_{i}=1-u_{i}^{1/\lambda_{i}},$ $i=1,2$,
we obtain%
\begin{align*}
\frac{1-\tau}{4\lambda_{1}\lambda_{2}}  &  =%
{\displaystyle\iint_{\left[  0,1\right]  ^{2}}}
\left(  1-r_{12}t_{1}t_{2}\right)  ^{-2\lambda}\left(  1-t_{1}\right)
^{2\lambda_{1}-1}\left(  1-t_{2}\right)  ^{2\lambda_{2}-1}\\
&  -\frac{\lambda}{\lambda_{1}}r_{12}\left(  1-r_{12}t_{1}t_{2}\right)
^{-2\lambda-1}\left(  1-t_{1}\right)  ^{2\lambda_{1}}t_{2}\left(
1-t_{2}\right)  ^{2\lambda_{2}-1}\\
&  -\frac{\lambda}{\lambda_{2}}r_{12}\left(  1-r_{12}t_{1}t_{2}\right)
^{-2\lambda-1}t_{1}\left(  1-t_{1}\right)  ^{2\lambda_{1}-1}\left(
1-t_{2}\right)  ^{2\lambda_{2}}\\
&  +\frac{\lambda^{2}}{\lambda_{1}\lambda_{2}}r_{12}^{2}\left(  1-r_{12}%
t_{1}t_{2}\right)  ^{-2\lambda-2}t_{1}\left(  1-t_{1}\right)  ^{2\lambda_{1}%
}t_{2}\left(  1-t_{2}\right)  ^{2\lambda_{2}}\,\text{d}t_{1}\text{d}t_{2}.
\end{align*}
By using equality $\left(  1-r_{12}t_{1}t_{2}\right)  ^{-2\lambda}=\sum
_{k=0}^{\infty}\left(  2\lambda\right)  _{k}\frac{\left(  r_{12}^{k}t_{1}%
^{k}t_{2}^{k}\right)  }{k!}$ in the last result, we obtain,%
\begin{align*}
\frac{1-\tau}{4\lambda_{1}\lambda_{2}}  &  =\sum_{k=0}^{\infty}\left(
2\lambda\right)  _{k}\frac{r_{12}^{k}}{k!}B\left(  k+1,2\lambda_{1}\right)
B\left(  k+1,2\lambda_{2}\right) \\
&  -\frac{\lambda}{\lambda_{1}}r_{12}\sum_{k=0}^{\infty}\left(  2\lambda
+1\right)  _{k}\frac{r_{12}^{k}}{k!}B\left(  k+1,2\lambda_{1}+1\right)
B\left(  k+2,2\lambda_{2}\right) \\
&  -\frac{\lambda}{\lambda_{2}}r_{12}\sum_{k=0}^{\infty}\left(  2\lambda
+1\right)  _{k}\frac{r_{12}^{k}}{k!}B\left(  k+2,2\lambda_{1}\right)  B\left(
k+1,2\lambda_{2}+1\right) \\
&  +\frac{\lambda^{2}}{\lambda_{1}\lambda_{2}}r_{12}^{2}\sum_{k=0}^{\infty
}\left(  2\lambda+2\right)  _{k}\frac{r_{12}^{k}}{k!}B\left(  k+2,2\lambda
_{1}+1\right)  B\left(  k+2,2\lambda_{2}+1\right)  .
\end{align*}
With the equality $B\left(  \alpha,\beta\right)  =\Gamma\left(  \alpha\right)
\Gamma\left(  \beta\right)  /\Gamma\left(  \alpha+\beta\right)  $ for
$\alpha,\beta>0,$ we obtain%
\begin{align*}
1-\tau &  =\sum_{k=0}^{\infty}\frac{\left(  1\right)  _{k}\left(  1\right)
_{k}\left(  2\lambda\right)  _{k}}{\left(  2\lambda_{1}+1\right)  _{k}\left(
2\lambda_{2}+1\right)  _{k}}\frac{r_{12}^{k}}{k!}\\
&  -\frac{4\lambda}{\left(  2\lambda_{1}+1\right)  \left(  2\lambda
_{2}+1\right)  }r_{12}\sum_{k=0}^{\infty}\frac{\left(  1\right)  _{k}\left(
2\right)  _{k}\left(  2\lambda+1\right)  _{k}}{\left(  2\lambda_{1}+2\right)
_{k}\left(  2\lambda_{2}+2\right)  _{k}}\frac{r_{12}^{k}}{k!}\\
&  +\frac{\lambda^{2}}{\left(  2\lambda_{1}+1\right)  \left(  2\lambda
_{2}+1\right)  \left(  \lambda_{1}+1\right)  \left(  \lambda_{2}+1\right)
}r_{12}^{2}\sum_{k=0}^{\infty}\frac{\left(  2\right)  _{k}\left(  2\right)
_{k}\left(  2\lambda+2\right)  _{k}}{\left(  2\lambda_{1}+3\right)
_{k}\left(  2\lambda_{2}+3\right)  _{k}}\frac{r_{12}^{k}}{k!}.
\end{align*}
Finally, by using equality (\ref{Hypergeomconf}), we have proved equality
(\ref{KtauHyper}).

\subsubsection{Spearman's rho}

By injecting the given copula in formula (\ref{bmfgLc}) in formula
((\ref{rho formula}), we obtain%

\begin{equation}
\varrho_{S}=12\int_{0}^{1}\int_{0}^{1}u_{1}u_{2}[1-r_{12}(1-u_{1}^{\frac
{1}{\lambda_{1}}})(1-u_{2}^{\frac{1}{\lambda_{2}}})]^{-\lambda}\,\text{d}%
u_{1}\text{d}u_{2}-3. \label{rho inter}%
\end{equation}
By changing the variables $v_{i}=u_{i}^{\frac{1}{\lambda_{i}}},$ $i=1,2$, we
obtain%
\begin{align}
&  \hspace{-3cm}\int_{0}^{1}\int_{0}^{1}u_{1}u_{2}[1-r_{12}\left(
(1-u_{1}^{\frac{1}{\lambda_{1}}})\right)  (1-u_{2}^{\frac{1}{\lambda_{2}}%
})]^{-\lambda}\,\text{d}u_{1}\text{d}u_{2}\nonumber\\
&  =\int_{0}^{1}\int_{0}^{1}\lambda_{1}\lambda_{2}v_{1}^{2\lambda_{1}-1}%
v_{2}^{2\lambda_{2}-1}\left[  1-r_{12}\left(  1-v_{1}\right)  \left(
1-v_{2}\right)  \right]  ^{-\lambda}\,\text{d}v_{1}\text{d}v_{2}\nonumber\\
&  =\lambda_{1}\lambda_{2}\sum_{k=0}^{\infty}\frac{\left(  \lambda\right)
_{k}}{k!}r_{12}^{k}\int_{0}^{1}v_{1}^{2\lambda_{1}-1}\left(  1-v_{1}\right)
^{k}\text{d}v_{1}\int_{0}^{1}v_{2}^{2\lambda_{2}-1}\left(  1-v_{2}\right)
^{k}\,\text{d}v_{2}\nonumber\\
&  =\lambda_{1}\lambda_{2}\sum_{k=0}^{\infty}\frac{\left(  \lambda\right)
_{k}}{k!}r_{12}^{k}B\left(  2\lambda_{1},k+1\right)  B\left(  2\lambda
_{2},k+1\right) \nonumber\\
&  =\frac{1}{4}\sum_{k=0}^{\infty}\frac{\left(  \lambda\right)  _{k}\left(
1\right)  _{k}^{2}}{\left(  2\lambda_{1}+1\right)  _{k}\left(  2\lambda
_{2}+1\right)  _{k}}\frac{r_{12}^{k}}{k!}\nonumber\\
&  =\frac{1}{4}F_{2}^{3}\left(  1,1,\lambda;2\lambda_{1}+1,2\lambda
_{2}+1;r_{12}\right)  \label{integral rho}%
\end{align}
By injecting equality (\ref{integral rho}) in equality (\ref{rho inter}), we
obtain equality (\ref{SrhoHyper1}). We remark that we have%
\[
F_{2}^{3}\left(  \lambda,1,1;2\lambda_{1}+1,2\lambda_{2}+1;r_{12}\right)
-1=\frac{\lambda}{\left(  2\lambda_{1}+1\right)  \left(  2\lambda
_{2}+1\right)  }r_{12}\left[  \sum_{k=0}^{\infty}\left(  k+1\right)
!\frac{\left(  \lambda+1\right)  _{k}}{\left(  2\lambda_{1}+2\right)
_{k}\left(  2\lambda_{2}+2\right)  _{k}}r_{12}^{k}\right]  .
\]
This gives equality (\ref{SrhoHyper2}).

\section{Acknowledgement}

I thank G\'{e}rard Letac for many helpful conversations.

\end{document}